\documentclass[11pt, final]{amsart}
\usepackage{amssymb}
\usepackage[mathscr]{eucal}

\theoremstyle{plain}
\newtheorem{theorem}{Theorem}
\newtheorem*{wtheorem}{Wigner's theorem}

\theoremstyle{remark}
\newtheorem*{remark}{Remark}

\newcommand\C{\mathbb C}
\newcommand\Mn{M_n(\mathbb F)}

\newcommand\Sn{S_n(\mathbb F)}
\newcommand\tr{\operatorname{tr}}

\begin{document}
\title[A Wigner-type theorem in Banach spaces]{A Wigner-type theorem on
symmetry transformations in Banach spaces}
\author{LAJOS MOLN\'AR}
\address{Institute of Mathematics and Informatics\\
         University of Debrecen\\
         4010 Debrecen, P.O.Box 12, Hungary}
\email{molnarl@math.klte.hu}
\thanks{  This research was supported from the following sources:\\
          1) Hungarian National Foundation for Scientific Research
          (OTKA), Grant No. T--030082, T--031995,\\
          2) A grant from the Ministry of Education, Hungary, Reg.
          No. FKFP 0349/2000}
\date{\today}
\begin{abstract}
We obtain an analogue of Wigner's classical theorem on symmetries
for Banach spaces. The proof is based on a result from the
theory of linear preservers. Moreover,
we present two other Wigner-type results for finite
dimensional linear spaces over general fields.
\end{abstract}
\maketitle

Wigner's theorem on symmetry transformations (sometimes called
unitary-antiunitary theorem) plays fundamental role in quantum
mechanics. This result can be formulated in several ways. For example,
in \cite[Theorem 3.1]{CVLL} the statement reads as follows.
In the sequel $P_1(H)$ denotes the set of all rank-one
(orthogonal) projections
(or, in the language of quantum mechanics, the set of all pure states)
on the Hilbert space $H$. We let $\tr$ stand for the usual
trace-functional.

\begin{wtheorem}
Let $H$ be a complex Hilbert space and let $\phi:P_1(H) \to P_1(H)$
be a bijective function for which
\begin{equation}\label{E1:harmadik}
\tr \phi(P)\phi(Q)=\tr PQ \qquad (P,Q \in P_1(H)).
\end{equation}
Then there exists an either unitary or antiunitary operator $U$ on $H$
such that $\phi$ is of the form
\[
\phi(P)=UPU^* \qquad (P\in P_1(H)).
\]
\end{wtheorem}

This formulation of Wigner's theorem makes us possible to formulate an
analogous theorem in the more general setting of Banach spaces, which we
state below.
Other Wigner-type theorems for Hilbert modules over matrix algebras
or for indefinite inner product spaces or for type II factors can be
found
in our recent papers \cite{MolJMP}, \cite{MolCMP},
\cite{MolIJTP}, respectively.

If $X$ is a Banach space, then $X'$ denotes the (topological) dual
of $X$. The set of all rank-one idempotents on $X$ (which are the
natural Banach space analogues of the projections) is denoted by
$I_1(X)$.

Now, our first result reads as follows.

\begin{theorem}\label{T:wigban1}
Let $X$ be a (real or complex) Banach space. Let $\phi:I_1(X) \to
I_1(X)$ be a bijective function for which
\begin{equation}\label{E1:miska}
\tr \phi(P)\phi(Q)=\tr PQ \qquad (P,Q \in I_1(X)).
\end{equation}
Then either there exists
a bijective bounded linear operator $A:X\to X$ such that
\[
\phi(P) = APA^{-1} \qquad (P\in I_1(X))
\]
or there exists a bijective bounded linear operator $C:X'\to X$ such
that
\[
\phi(P) = CP'C^{-1} \qquad (P\in I_1(X).
\]
\end{theorem}

Before the proof we need some additional notation and definitions.
If $X$ is a Banach space, then
$B(X)$ stands for the algebra of all bounded linear operators
on $X$ and $F(X)$ denotes the ideal of all finite rank elements in
$B(X)$. The Banach space adjoint of the operator $A\in B(X)$ is
denoted by $A'$. If $x\in X$, $f\in X'$, then $x\otimes f$ is the
operator defined by
\[
(x\otimes f)(y)=f(y)x \qquad (y\in X).
\]
It is easy to see that $x\otimes f\in I_1(X)$ if and
only if $f(x)=1$.
Clearly, every operator of the form
$\sum_{i=1}^n x_i \otimes f_i$ belongs to $F(X)$, and,
conversely, every finite
rank operator $A\in F(X)$ can be written in the form
\[
A=\sum_{i=1}^n x_i \otimes f_i
\]
with some $x_1, \ldots, x_n \in X$ and $f_1, \ldots, f_n \in X'$.
On the elements of $F(X)$ we define the trace-functional by
\[
\tr \sum_{i=1}^n x_i \otimes f_i=
\sum_{i=1}^n f_i(x_i).
\]
One can see that this functional is well-defined and,
in the Hilbert space case, it gives us the usual trace. It
is not hard to see that $\tr$ is a linear functional on $F(X)$ with the
property that
\[
\tr TA=\tr AT
\]
for every $A\in F(X)$ and $T\in B(X)$ (cf. \cite[Section B.1]{Pie}).
In particular, this shows that any map $\phi$ of any of the forms
appearing in our statement above necessarily satisfies \eqref{E1:miska}.

\begin{proof}[Proof of Theorem~\ref{T:wigban1}]
The idea of the proof is very simple. We first extend $\phi$ from
$I_1(X)$ to a linear transformation on the operator algebra $F(X)$ and
then apply a result of Omladi\v c and \v Semrl on linear preservers.
Notice that our approach is completely different from the usual proofs
of Wigner's theorem and, in particular, this is the case also with the
proof presented in \cite{CVLL}.

First suppose that $X$ is a complex Banach space.
We define
\[
\Phi(\sum_{i=1}^n \lambda_i P_i)=
\sum_{i=1}^n \lambda_i \phi(P_i)
\]
whenever $P_1, \ldots , P_n$ are rank-one idempotents and $\lambda_1,
\ldots ,\lambda_n$ are scalars. We assert that $\Phi$ is well-defined.
To see this, we prove that
\begin{equation}\label{E1:elso}
\sum_{i=1}^n \lambda_i P_i=0
\end{equation}
implies
\[
\sum_{i=1}^n \lambda_i \phi(P_i)=0.
\]
Indeed, from \eqref{E1:elso} we obtain that
\[
\sum_{i=1}^n \lambda_i P_iQ=0 \qquad (Q\in I_1(X)).
\]
By the linearity of the trace and the equality \eqref{E1:miska}, it
follows that
\[
\tr(\sum_{i=1}^n \lambda_i \phi(P_i)\phi(Q))=0 \qquad (Q\in I_1(X)).
\]
Since $\phi$ maps onto $I_1(X)$, we can infer that
\[
\sum_{i=1}^n \lambda_i \phi(P_i)=0.
\]
So $\Phi$ is well-defined.
Clearly, $\Phi$ is a linear transformation on $F(X)$. Since every
matrix $A\in M_n(\C)$ is a linear combination of rank-one idempotents,
it follows that every finite-rank operator belongs to the linear span of
$I_1(X)$.
This yields that $\Phi$ is defined on the whole $F(X)$ and
maps onto $F(X)$ (in fact, one can prove that
it is injective as well).

So we have a surjective linear transformation $\Phi$ on $F(X)$ which
preserves the rank-one idempotents. Now, we can apply a result of
Omladi\v c and \v Semrl describing the form of all such maps.
In view of \cite[Main Theorem]{OmlSem}, we distinguish two cases. First
suppose that $X$ is infinite-dimensional. Then for the form of $\Phi$,
and hence for the form of $\phi$, we have the following
two possibilities:
\begin{itemize}
\item[(i)]
There exists a bijective bounded linear operator $A:X\to X$ such that
\[
\phi(P) = APA^{-1} \qquad (P\in I_1(X)).
\]
\item[(ii)]
There exists a bijective bounded linear operator $C:X'\to X$ such that
\[
\phi(P) = CP'C^{-1} \qquad (P\in I_1(X)).
\]
\end{itemize}
Recall that our map $\Phi$ is linear. This is the reason that the two
remaining possibilities in \cite[Main Theorem]{OmlSem} do not appear
here.
As for the finite dimensional case, the form of all surjective linear
maps
on $M_n(\C)$ preserving rank-one idempotents is given in \cite[Theorem
4.5]{OmlSem}. Since every linear transformation on $M_n(\C)$ is
continuous and the only continuous ring automorphisms of $\C$ are the
identity and the conjugation, we obtain from that result that in this
case our map $\phi$ is either of the form
$P\mapsto APA^{-1}$ or of the form
$P\mapsto AP^{t}A^{-1}$ (${}^t$ stands for the transpose) with some
nonsingular matrix $A\in M_n(\C)$.
But under the identification of matrices and operators, we obtain the
forms appearing in (i) and (ii) once again.

If $X$ is a real Banach space, then one can argue in the same way but
since
in that case \cite[Main Theorem]{OmlSem} holds without any assumption on
the dimension of $X$, there is no need to refer to \cite[Theorem
4.5]{OmlSem}.
\end{proof}

\begin{remark}
Using our quite algebraic approach presented above
we can now give a short proof of Wigner's original theorem.
Let $H$ be a complex Hilbert space and let $\phi: P_1(H) \to P_1(H)$ be
a bijective function which preserves the transition probabilities, that
is, satisfies \eqref{E1:harmadik}.
One can extend $\phi$ to a linear transformation $\Phi$ on
$F(H)$ in a similar way as in the proof of our previous theorem. Just as
there,
we find that $\Phi$ is bijective. On the other hand, for any $P,Q\in
P_1(H)$ we have $PQ=QP=0$ if and only if $\tr PQ=0$.
Therefore, $\Phi$ preserves the orthogonality between rank-one
projections which, by the linearity of $\Phi$, implies that $\Phi$
sends projections to projections. It is now an easy algebraic argument
to verify that $\Phi$ is a Jordan *-automorphism of $F(H)$, that is,
$\Phi$ satisfies $\Phi(A^2)=\Phi(A)^2$ and $\Phi(A^*)=\Phi(A)^*$ for
every $A\in F(H)$ (see \cite[Remark 2.2]{BreSem}).
By a classical result of Herstein \cite{Her}, $\Phi$ is either a
*-automorphism or a *-antiautomorphism of $F(H)$. But the forms
of those morphisms are well-known
(see, for example, \cite[Proposition]{MolJNG}). Namely, we have an
either unitary or antiunitary operator $U$ on $H$ such that
\[
\phi(P)=UPU^* \qquad (P\in P_1(H)).
\]
This completes the proof of Wigner's theorem.
\end{remark}

In the finite dimensional case we can get rid of the assumption on the
bijectivity of the transformation $\phi$. In fact, we have the following
Wigner-type result for matrix algebras over general fields.

\begin{theorem}\label{T:wigban2}
Let $\mathbb F$ be a field of characteristic different from 2 and
let $n\in \mathbb N$.
Suppose $\phi$ is a transformation on the set $I_1(\mathbb F^n)$ of all
rank-one idempotents in $\Mn$ into itself with the property that
\[
\tr \phi(P)\phi(Q)=\tr PQ \qquad (P,Q \in I_1(\mathbb F^n)).
\]
Then there exists a nonsingular matrix $A\in \Mn$
such that $\phi$ is either of the form
\[
\phi(P)=APA^{-1} \qquad
(P\in I_1(\mathbb F^n))
\]
or of the form
\[
\phi(P)=AP^{t} A^{-1} \qquad
(P\in I_1(\mathbb F^n)).
\]
\end{theorem}

\begin{proof}
We first show that the range of $\phi$ linearly generates $\Mn$.
As usual, denote by $E_{ij}\in \Mn$
the matrix whose $ij$ entry is 1 and its all other entries are 0.
Define
\[
E_{ij}'=
\begin{cases}
\phi(E_{ii}+E_{ij})-\phi(E_{ii}) &\text{if $i\neq j$;}\\
\phi(E_{ii})                     &\text{if $i= j$.}
\end{cases}
\]
Suppose that
\[
\sum_{i,j} \lambda_{ij}E_{ij}' =0
\]
for some scalars $\lambda_{ij}\in \mathbb F$.
Fix indices $k,l\in \{ 1, \ldots, n\}$. We have
\[
\sum_{i,j} \lambda_{ij}E_{ij}'E_{kl}' =0.
\]
Taking trace, we obtain
\[
\sum_{i,j} \lambda_{ij}\tr E_{ij}'E_{kl}' =0.
\]
By the preserving property of $\phi$, it follows that
\[
\sum_{i,j} \lambda_{ij}\tr E_{ij}E_{kl} =0.
\]
Since $E_{ij}E_{kl}=\delta_{jk} E_{il}$, from this equality we easily
deduce that $\lambda_{lk}=0$. As $k,l$ were arbitrary, it follows that
the matrices $E_{ij}'$, $i,j\in \{ 1, \ldots, n\}$ form a linearly
independent set in
$\Mn$. This implies that the range of $\phi$ linearly generates $\Mn$.

Similarly to the proof of Theorem~\ref{T:wigban1}
we define a transformation $\Phi :\Mn \to \Mn$ by
\[
\Phi(\sum_i \lambda_i P_i)=\sum_i \lambda_i \phi(P_i),
\]
where the $P_i$'s are rank-one idempotents and the $\lambda_i$'s are
scalars. First, $\Phi$ is well-defined. Indeed, let
$Q_j$ be rank-one idempotents and $\mu_j\in \mathbb F$ be such that
$\sum_i \lambda_i P_i=\sum_j \mu_j Q_j$. Then we have
\[
\sum_i \lambda_i \tr P_iE_{kl}=
\tr (\sum_i \lambda_i P_i)E_{kl}=
\tr (\sum_j \mu_j Q_j)E_{kl}=
\sum_j \mu_j \tr Q_jE_{kl}.
\]
By the preserving property of $\phi$ we infer that
\[
\sum_i \lambda_i \tr \phi(P_i)E_{kl}'=
\sum_j \mu_j    \tr \phi(Q_j)E_{kl}'.
\]
This implies that
\[
\tr (\sum_i \lambda_i \phi(P_i))E_{kl}'=
\tr (\sum_j \mu_j \phi(Q_j))E_{kl}'
\]
and, as the matrices $E_{kl}'$ linearly generate $\Mn$, we can conclude
that
$\sum_i \lambda_i \phi(P_i) =\sum_j \mu_j \phi(Q_j)$. Therefore, $\Phi$
is well-defined. Since the rank-one idempotents linearly generate
$\Mn$, we obtain that
$\Phi$ is a linear transformation from $\Mn$ into itself which preseves
the rank-one idempotents. The form of such transformations is described
in \cite[Theorem 3]{CL} and this gives us the form of $\phi$.
\end{proof}

Wigner's classical result concerns transformations on the set of all
rank-one (orthogonal) projections
on a Hilbert space. In case the Hilbert space in question is finite
dimensional and real, then the result reduces to the description
of all bijective transformations on the set of all symmetric rank-one
idempotents in $M_n (\mathbb R)$ which preserve the trace of the
product. In our last assertion using a result on linear
preservers once again, we can generalize this statement for the case of
general fields.

\begin{theorem}
Let $\mathbb F$ be an algebraically closed field of characteristic
different from 2 and let $n\in \mathbb N$.
Suppose $\phi$ is a transformation on the set $P_1(\mathbb F^n)$ of all
rank-one idempotents
in $\Sn$ (the set of all symmetric elements of $\Mn$) into itself
with the property that
\[
\tr \phi(P)\phi(Q)=\tr PQ \qquad (P,Q \in P_1(\mathbb F^n)).
\]
Then there exists an orthogonal matrix $U\in \Mn$
such that $\phi$ is either of the form
\[
\phi(P)=UPU^{-1} \qquad
(P\in P_1(\mathbb F^n))
\]
or of the form
\[
\phi(P)=UP^{t} U^{-1} \qquad
(P\in P_1(\mathbb F^n)).
\]
\end{theorem}

\begin{proof}
The proof is quite similar to the proof of our previous result.
Define
\[
E_{ij}'=
\begin{cases}
2\phi((E_{ii}+E_{jj}+E_{ij}+E_{ji})/2)-(\phi(E_{ii})+\phi(E_{jj}))
                                 &\text{if $i< j$;}\\
\phi(E_{ii})                     &\text{if $i= j$.}
\end{cases}
\]
Just as in the proof of Theorem~\ref{T:wigban2}, one can show that
the matrices $E_{ij}'$, $i\leq j$ are linearly independent in $\Sn$.
This gives us that the range of $\phi$ linearly generates $\Sn$.

Next we define
\[
\Phi(\sum_i \lambda_i P_i)= \sum_i \lambda_i \phi(P_i)
\]
for every finite system $\lambda_i$ of scalars and symmetric rank-one
idempotents $P_i$. One can prove that $\Phi$ is well-defined in a
way very similar to the corresponding part of
the proof of Theorem~\ref{T:wigban2}.
Since the symmetric rank-one idempotents linearly generate $\Sn$,
$\Phi$ is a linear transformation from $\Sn$ into
itself which preserves the rank-one idempotents in $\Sn$. The form of
such transformations on $\Sn$ is described in \cite[Theorem 4]{CL}
and this gives us the form of $\phi$.
\end{proof}

\end{document}